\documentstyle[11pt,twoside]{article}
\include{graphicx}
\title{The asymptotics of the Struve function ${\bf H}_\nu(z)$ for large complex order and argument}
\author{\sc R. B.\ Paris \\
{\em Division of Computing and Mathematics}, \\
{\em University of Abertay Dundee, Dundee DD1 1HG, UK}
}
\begin{document}
\def\f#1#2{\mbox{${\textstyle \frac{#1}{#2}}$}}
\def\dfrac#1#2{\displaystyle{\frac{#1}{#2}}}
\def\boldal{\mbox{\boldmath $\alpha$}}
{\newcommand{\Sgoth}{S\;\!\!\!\!\!/}
\newcommand{\bee}{\begin{equation}}
\newcommand{\ee}{\end{equation}}
\newcommand{\lam}{\lambda}
\newcommand{\ka}{\kappa}
\newcommand{\al}{\alpha}
\newcommand{\fr}{\frac{1}{2}}
\newcommand{\fs}{\f{1}{2}}
\newcommand{\g}{\Gamma}
\newcommand{\br}{\biggr}
\newcommand{\bl}{\biggl}
\newcommand{\ra}{\rightarrow}
\newcommand{\mbint}{\frac{1}{2\pi i}\int_{c-\infty i}^{c+\infty i}}
\newcommand{\mbcint}{\frac{1}{2\pi i}\int_C}
\newcommand{\mboint}{\frac{1}{2\pi i}\int_{-\infty i}^{\infty i}}
\newcommand{\gtwid}{\raisebox{-.8ex}{\mbox{$\stackrel{\textstyle >}{\sim}$}}}
\newcommand{\ltwid}{\raisebox{-.8ex}{\mbox{$\stackrel{\textstyle <}{\sim}$}}}
\renewcommand{\topfraction}{0.9}
\renewcommand{\bottomfraction}{0.9}
\renewcommand{\textfraction}{0.05}
\newcommand{\mcol}{\multicolumn}
\date{}
\maketitle
\pagestyle{myheadings}
\markboth{\hfill \sc R. B.\ Paris  \hfill}
{\hfill \sc  Asymptotics of the Struve function\hfill}
\begin{abstract}
We re-examine the asymptotic expansion of the Struve function ${\bf H}_\nu(z)$ for large complex values of $\nu$ and $z$ satisfying $|\arg\,\nu|\leq\fs\pi$ and $|\arg\,z|<\fs\pi$. Watson's
analysis \cite[\S 10.43]{W} covers only the case of $\nu$ and $z$ of the same phase with $\nu/z$ in the intervals $(0,1)$ and $(1,\infty)$. The domains in the complex $\nu/z$-plane where the expansion takes on different forms are obtained.
\vspace{0.4cm}

\vspace{0.3cm}

\end{abstract}

\vspace{0.3cm}

\noindent $\,$\hrulefill $\,$

\vspace{0.2cm}

\begin{center}
{\bf 1. \  Introduction}
\end{center}
\setcounter{section}{1}
\setcounter{equation}{0}
\renewcommand{\theequation}{\arabic{section}.\arabic{equation}}
The Struve function ${\bf H}_\nu(z)$ is a particular solution of the inhomogeneous Bessel equation
\[\frac{d^2w(z)}{dz^2}+\frac{1}{z}\,\frac{dw(z)}{dz}+\bl(1-\frac{\nu^2}{z^2}\br)\,w(z)=\frac{(\fs z)^{\nu-1}}{\surd\pi\,\g(\nu+\fs)}\]
which possesses the series expansion
\bee\label{e11}
{\bf H}_\nu(z)=(\fs z)^{\nu+1} \sum_{n=0}^\infty \frac{(-)^n (\fs z)^{2n}}{\g(n+\f{3}{2}) \g(n+\nu+\f{3}{2})}
\ee
valid for all finite $z$.

An integral representation, valid when $\Re (\nu)>-\fs$, is given by \cite[p.~330]{W} as 
\[J_\nu(z)\pm i {\bf H}_\nu(z)=\frac{2(\fs z)^\nu}{\surd\pi\,\g(\nu+\fs)} \int_0^1 e^{\pm izt}(1-t^2)^{\nu-\fr} \,dt,\]
where $J_\nu(z)$ is the usual Bessel function. Upon replacement of the variable $t$ by $\pm iu$, we obtain
\bee\label{e12}
{\bf H}_\nu(z)\pm iJ_\nu(z)=\frac{2(\fs z)^\nu}{\surd\pi\,\g(\nu+\fs)} \int_0^{\pm i} e^{-zu}(1+u^2)^{\nu-\fr}\, du \quad (\Re (\nu)>-\fs).
\ee
The integration path corresponding to the upper sign in (\ref{e12}) can be deformed to pass along the positive real axis to $+\infty$ and back to the point $i$ along the parallel path $i+u$ ($0\leq u\leq\infty$). The contribution from the path $(i+\infty, i]$ is equal to $iH_\nu^{(2)}(z)$, where $H_\nu$ is the Hankel function; see \cite[p.~166]{W}. Thus we find the alternative representation \cite[p.~292]{DLMF}
\bee\label{e13}
{\bf H}_\nu(z)-Y_\nu(z)=\frac{2 (\fs z)^\nu}{\surd\pi\,\g(\nu+\fs)} \int_0^\infty e^{-zu}(1+u^2)^{\nu-\fr}\,du
\ee
valid\footnote{Suitable rotation of the integration path through an acute angle enables the validity of (\ref{e13}) to be extended to the wider sector $|\arg\,z|<\pi$; see \cite[p.~331]{W}.} for unrestricted $\nu$ and $|\arg\,z|<\fs\pi$, where $Y_\nu(z)$ denotes the Bessel function of the second kind.

Here we shall consider the asymptotic expansion of ${\bf H}_\nu(z)$ for large complex values of $\nu$ and $z$ satisfying $|\arg\,\nu|\leq\fs\pi$ and $|\arg\,z|<\fs\pi$. 
Values of $\arg\,z$ outside this range can be dealt with by means of the continuation formula
\[{\bf H}_\nu(ze^{\pi mi})=e^{\pi mi(\nu+1)}\,{\bf H}_\nu(z), \qquad m=\pm 1, \pm 2 \ldots \]
obtained from (\ref{e11}).
\vspace{0.6cm}

\begin{center}
{\bf 2. \ Asymptotic expansion when $z>0$}
\end{center}
\setcounter{section}{2}
\setcounter{equation}{0}
\renewcommand{\theequation}{\arabic{section}.\arabic{equation}}
We set 
\[q:=\nu/z=\alpha+i\beta,\qquad \theta:=\arg\,z,\qquad \omega:=\arg\,q.\]
In view of (\ref{e12}) and (\ref{e13}), we are led to the consideration of the integral
\bee\label{e21}
\int_C e^{-|z|\tau} \frac{du}{\sqrt{1+u^2}}\qquad \tau:=e^{i\theta} \{u-q \log (1+u^2)\},
\ee
where $C$ is a suitably chosen path in the $u$-plane. 

Saddle points are situated at $d\tau/du=0$; that is, at the points
\[u_\pm=q\pm\sqrt{q^2-1}.\]
We shall refer to these saddles as $S_1$ (uuper sign) and $S_2$ (lower sign). 
Inversion of (\ref{e21}) in the form $u=\sum_{k=1}^\infty a_k(\tau e^{-i\theta})^k$, where $a_0=1$, shows that
\[\frac{1}{\sqrt{1+u^2}}\,\frac{du}{d\tau}=e^{-i\theta}\sum_{k=0}^\infty c_k(q)(\tau e^{-i\theta})^k\]
valid in a disc centered at $\tau=0$ of radius determined by the nearest singularity corresponding to the saddles $S_1$ or $S_2$ (or both).
The values of the coefficients $c_k(q)$ ($0\leq k\leq 10$) are listed in Table 1; see also \cite[p.~203]{D}. 

Watson \cite[\S 10.43]{W} has considered the two cases (i) $q\in [1,\infty)$ and (ii) $q\in (0,1)$ when $z>0$ ($\theta=0$).
The steepest descent paths emanating from the origin in the complex $u$-plane in these two cases are shown in Fig.~1; branch cuts have been taken along the segments of the imaginary axis $[\pm i, \pm\infty i)$. In case (i), the desired path $C$ consists of the real axis between the origin and the saddle $S_2$ and then {\it either} along the arc\footnote{When $q=1$, the saddles $S_1$ and $S_2$ form a double saddle at $u=1$. In this case, the path $C$ consists of the real axis $0\leq u\leq 1$ followed by similar arcs to the points $u=\pm i$.} to the branch point at $u=i$ {\it or} along the arc to the branch point at $u=-i$. In case (ii), the path $C$ from the origin coincides with the positive real axis and passes to $+\infty$. In both cases $\tau$ increases monotonically from 0 to $+\infty$ as we traverse these paths.
\begin{figure}[t]
	\begin{center}
	{\tiny($a$)}\ \includegraphics[width=0.35\textwidth]{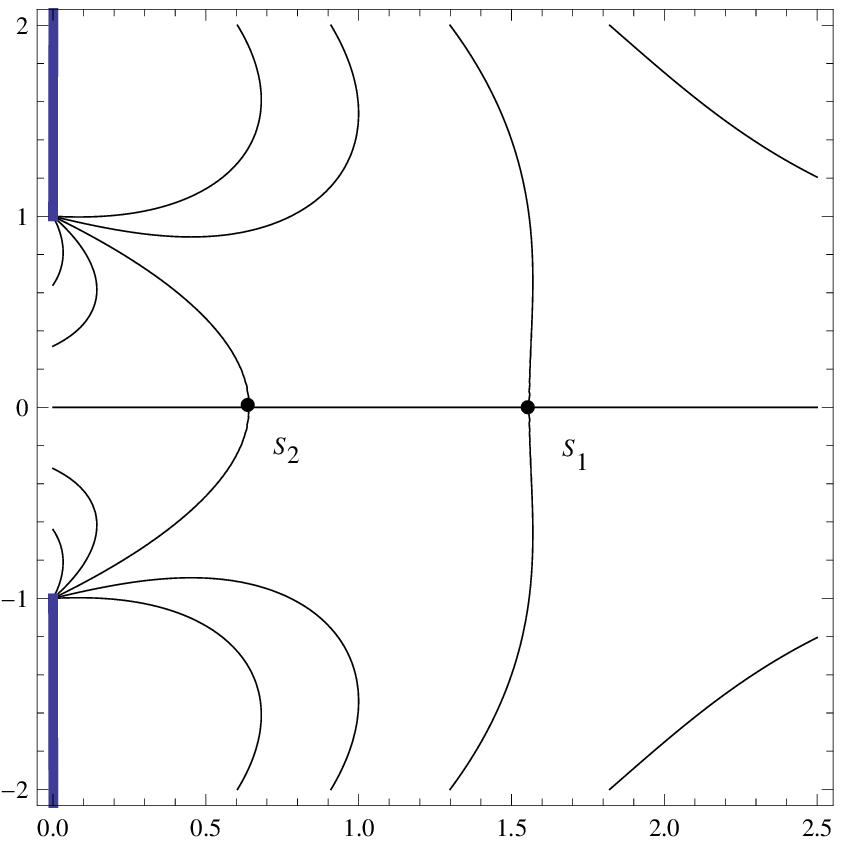}\qquad {\tiny($b$)} \includegraphics[width=0.35\textwidth]{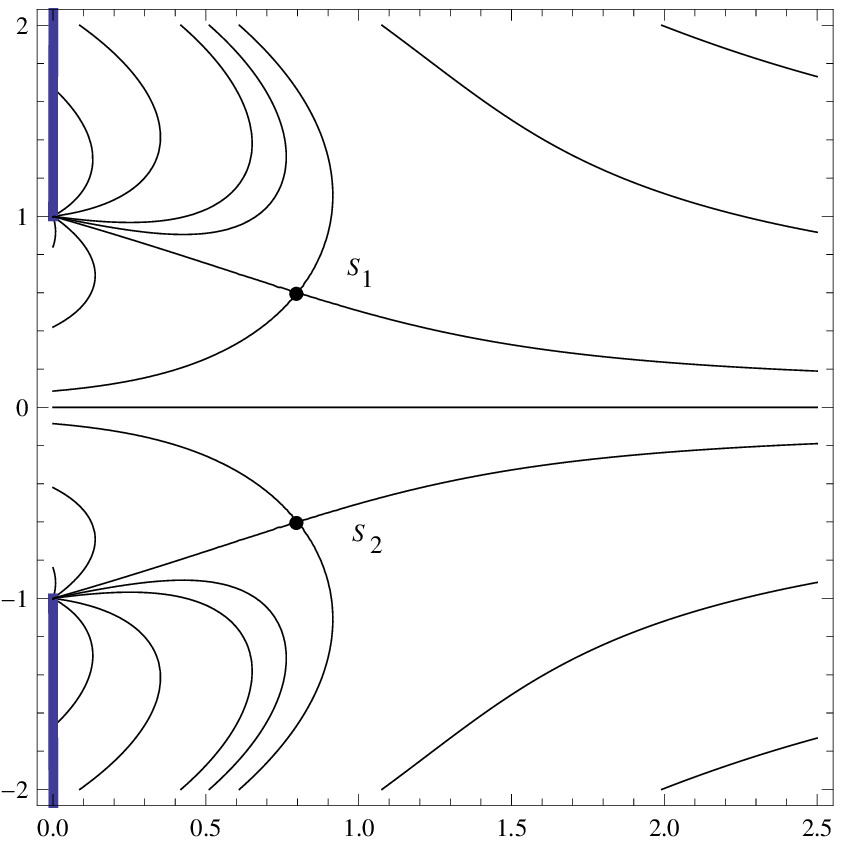}
	\caption{\small{The steepest paths when $\theta=0$: (a) when $q\in(1,\infty)$ and (b) when $q\in(0,1)$. The heavy dots indicate the saddle points and the heavy lines denote the branch cuts.}}
	\end{center}
\end{figure}

\begin{table}[h]
\caption{\footnotesize{The coefficients $c_k(q)$ for $0\leq k\leq 10$.}}
\begin{center}
\begin{tabular}{|l|l|}
\hline
&\\[-0.3cm]
\mcol{1}{|c|}{$k$}& \mcol{1}{c|}{$c_k(q)$}  \\
[.1cm]\hline
&\\[-0.3cm]
0 & $1$\\\
&\\[-0.4cm]
1 & $2q$\\
&\\[-0.4cm]
2 & $6q^2-\fs$\\
&\\[-0.4cm]
3 & $20q^3-4q$\\
&\\[-0.4cm]
4 & $70q^4-\f{45}{2}q^2+\f{3}{8}$\\
&\\[-0.4cm]
5 & $252q^5-112q^3+\f{23}{4}q$\\
&\\[-0.4cm]
6 & $924q^6-525q^4+\f{301}{6}q^2-\f{5}{16}$\\
&\\[-0.4cm]
7 & $3432q^7-2376q^5+345q^3-\f{22}{3}q$\\
&\\[-0.4cm]
8 & $12870q^8-\f{21021}{2}q^6+\f{16665}{8}q^4-\f{1425}{16}q^2+\f{35}{128}$\\
&\\[-0.4cm]
9 & $48620q^9-45760q^7+\f{139139}{12}q^5-\f{1595}{2}q^3+\f{563}{64}q$\\
&\\[-0.4cm]
10& $184756q^{10}-196911q^8+61061q^6-\f{287287}{48}q^4+\f{133529}{960}q^2-\f{63}{256}$\\
[.2cm]\hline
\end{tabular}
\end{center}
\end{table}

Then in case (i) we find 
\[\int_0^{\pm i} e^{-zu} (1+u^2)^{\nu-\fr}du=\int_0^\infty e^{-z\tau}\,\bl(\frac{1}{\sqrt{1+u^2}}\,\frac{du}{d\tau}\br)d\tau \sim \sum_{k=0}^\infty \frac{c_k(q) \g(k+1)}{z^{k+1}}\]
for $z\ra+\infty$.
Hence, for large real $\nu$ and $z$ with $\nu/z\in[1,\infty)$ (when the deformed path $C$ terminates at the branch points $u=\pm i$), we have from (\ref{e12})  
\bee\label{e22}
{\bf H}_\nu(z)\pm i J_\nu(z)\sim \frac{(\fs z)^{\nu-1}}{\sqrt{\pi} \g(\nu+\fs)}\,\sum_{k=0}^\infty \frac{c_k(q) \g(k+1)}{z^k}~,
\ee
respectively.
Similarly, for $\nu/z\in(0,1)$ (when the path $C$ passes to $+\infty$ along the real axis), we have from (\ref{e13})
\bee\label{e23}
{\bf H}_\nu(z)-Y_\nu(z)\sim \frac{(\fs z)^{\nu-1}}{\sqrt{\pi} \g(\nu+\fs)}\,\sum_{k=0}^\infty \frac{c_k(q) \g(k+1)}{z^k}~.
\ee
These are the results given in \cite[\S 10.43]{W}; see also the discussion in Section 3.

When $\nu$ is allowed to take on complex values with $z>0$, the steepest descent paths in Fig.~1 undergo a progressive change. Recalling that $q=\alpha+i\beta$, we find that as
$\beta$ increases from zero when $\alpha\in (0,1)$ the steepest descent path from the origin $\Im \tau=0$
becomes increasingly deformed in the upper-half plane, until at a critical value $\beta=\beta^*$ this path connects with the saddle $S_1$. For example, when $\alpha=0.80$ the critical value is $\beta^*\doteq 0.143900$.
Then, the path $\Im \tau=0$ passes to infinity when $\beta<\beta^*$, connects with $S_1$ when $\beta=\beta^*$ and approaches the branch point at $u=i$ (possibly spiralling onto different Riemann sheets) when $\beta>\beta^*$. 
An analogous transition occurs when $\beta<0$ at $\beta=-\beta^*$, with the saddle $S_1$ replaced by $S_2$.
When $\alpha>1$, the steepest path $\Im \tau=0$  passes to $u=i$ when $\beta>0$, and to $u=-i$ when $\beta<0$, without undergoing any transition as $\beta$ increases.

The transitions that occur when $z>0$ and $|\arg\,\nu|\leq\fs\pi$ are summarised in Fig.~2(a). This shows the three curves in the complex $q$-plane, on which a transition takes place, that emanate from the point $P$ (corresponding to $q=1$). The curves in the upper and lower half-planes are conjugate curves with the third being the segment $[1,\infty)$ of the real $q$-axis. In the domain numbered 1 (between the conjugate curves and the imaginary $q$-axis), the path $C$ passes to $\infty$ and the expansion (\ref{e23}) applies. In the domain numbered 2, the path $C$ terminates at $u=+i$ and the expansion (\ref{e22}) applies with the upper sign; in the domain numbered 3, the terminal point is $u=-i$ and the expansion (\ref{e22}) applies with the lower sign. For $q$ situated on these curves the transition is associated with a Stokes phenomenon; see below.
\vspace{0.6cm}

\begin{center}
{\bf 3. \ Asymptotic expansion for complex $z$}
\end{center}
\setcounter{section}{3}
\setcounter{equation}{0}
\renewcommand{\theequation}{\arabic{section}.\arabic{equation}}
When $z$ is complex ($\theta\neq 0$) the transition curves in the sector of the $q$-plane given by\footnote{This sector corresponds to $|\arg\,\nu|\leq\fs\pi$ and $|\arg\,z|<\fs\pi$.} $(-\fs\pi-\theta, \fs\pi-\theta)$ are $\theta$-dependent. In Fig.2(b)--(d) we show these curves for $\theta/\pi=0.10$, $0.20$ and $0.30$. The curves for $\theta<0$ are the conjugate of those for $\theta>0$. The point $P$ corresponds to the case when the steepest descent path from the origin connects with {\it both\/} saddles $S_1$ and $S_2$. The point labelled $Q$ is the intercept of the lower curve with the positive $q$-axis. Values of $q$ at $P$ and $Q$ are presented in Table 2 for different $\theta$.
\begin{table}[th]
\caption{\footnotesize{The coordinates of the triple point $P$ and the intercept $Q$ on the real $q$-axis as a function of $\theta$.}}
\begin{center}
\begin{tabular}{|l|ll||l|ll|}
\hline
&&&&&\\[-0.3cm]
\mcol{1}{|c|}{$\theta/\pi$}& \mcol{1}{c}{$P$} & \mcol{1}{c||}{$Q$} & \mcol{1}{c|}{$\theta/\pi$} & \mcol{1}{c}{$P$} & \mcol{1}{c|}{$Q$}\\
[.1cm]\hline
&&&&&\\[-0.3cm]
0    & $1$ & $1$                      & 0.30 & $1.08553+1.38238i$ & $0.27561$\\
0.05 & $0.96385+0.08606i$ & $0.83360$ & 0.35 & $1.36479+2.60425i$ & $0.18575$\\ 
0.10 & $0.93778+0.18745i$ & $0.70952$ & 0.40 & $2.36238+7.23955i$ & $0.10710$\\
0.20 & $0.93437+0.53249i$ & $0.48057$ & 0.42 & $3.72266+14.4826i$ & $0.07942$\\
0.25 & $0.97678+0.84047i$ & $0.37449$ & 0.45 & $16.4886+104.102i$ & $0.04275$\\
[.2cm]\hline
\end{tabular}
\end{center}
\end{table}

As in the case $\theta=0$ in Fig.~2(a), for $q$-values in domain 1 the endpoint of the steepest descent path from the origin terminates at infinity, whereas those situated in domains 2 and 3 pass to the branch points (possibly spiralling onto adjacent Riemann surfaces) at $u=\pm i$, respectively. As one crosses one of these curves, say from domain 1 to domain 2, there is a change in the endpoint 
via a Stokes phenomenon. Examples of the steepest descent paths when $\theta=0.10\pi$ on the three curves labelled $PA$, $PB$ and $PC$ in Fig.~2(b), and at $P$, are shown in Fig.~3 demonstrating that on each curve the change of endpoint is associated with a Stokes phenomenon.
\begin{figure}[t]
	\begin{center}
	{\tiny($a$)}\ \includegraphics[width=0.35\textwidth]{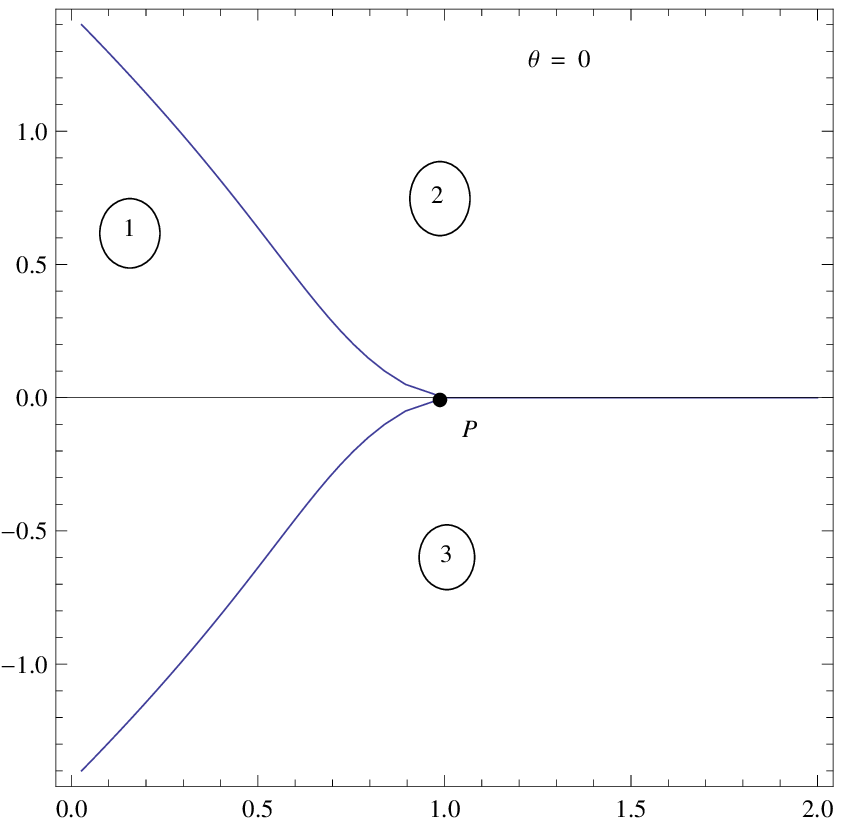}\qquad {\tiny($b$)} \includegraphics[width=0.35\textwidth]{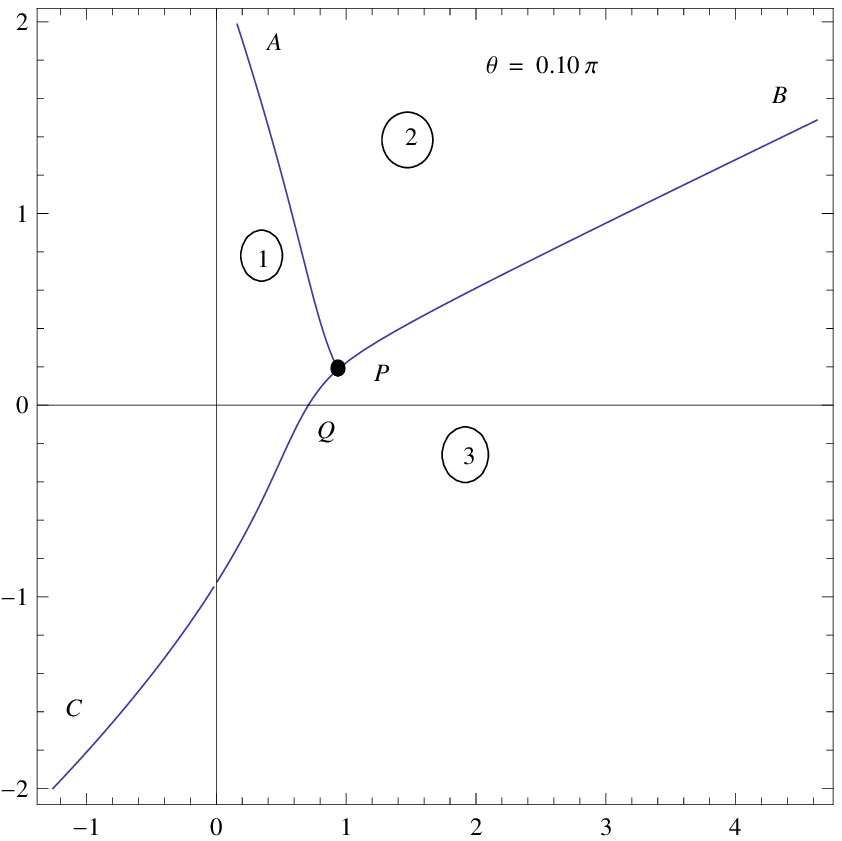} \vspace{.4cm}

{\tiny($c$)}\includegraphics[width=0.35\textwidth]{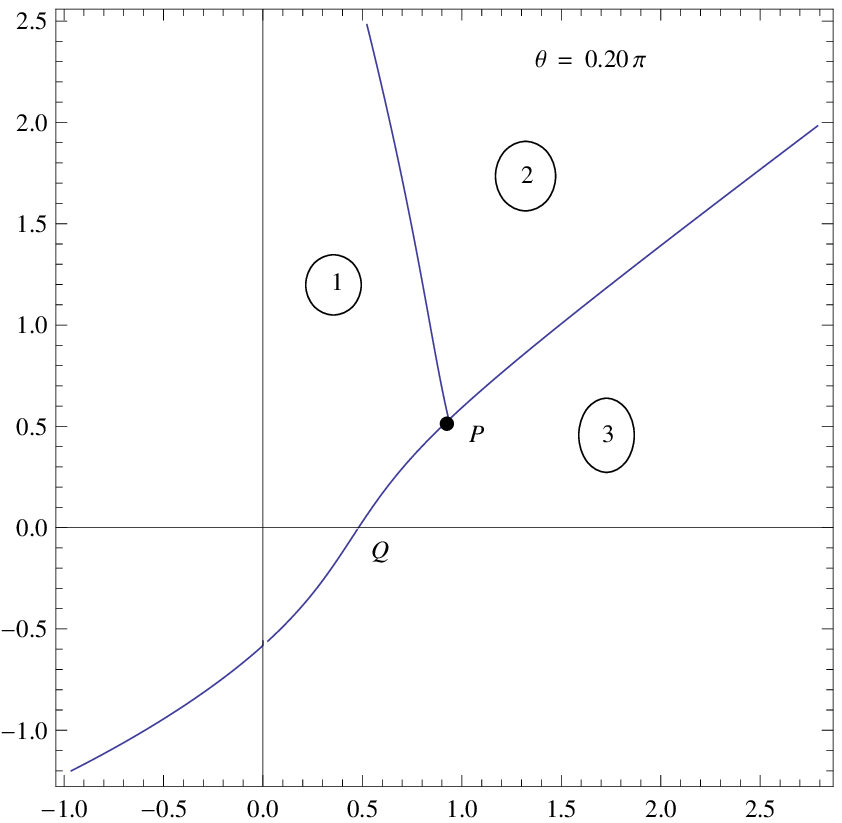}\qquad {\tiny($d$)}\includegraphics[width=0.35\textwidth]{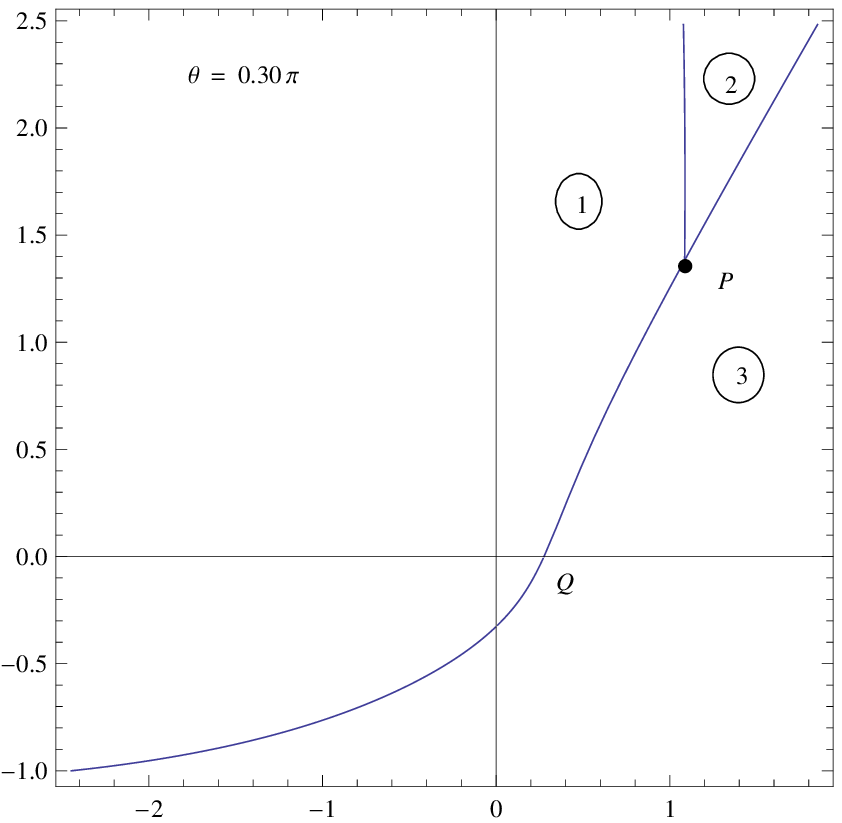}\\
\caption{\small{The domains in the sector of the $q$-plane bounded by $-\fs\pi-\theta<\omega<\fs\pi-\theta$ showing the termination points of the steepest descent path from the origin: (a) $\theta=0$, (b) $\theta=0.10\pi$, (c) $\theta=0.20\pi$ and (d) $\theta=0.30\pi$.	The termination point in domain 1 is at infinity and that in domains 2 and 3 is at $\pm i$, respectively.}}
\end{center}

\end{figure}
\begin{figure}[t]
	\begin{center}
	{\tiny($a$)}\ \includegraphics[width=0.28\textwidth]{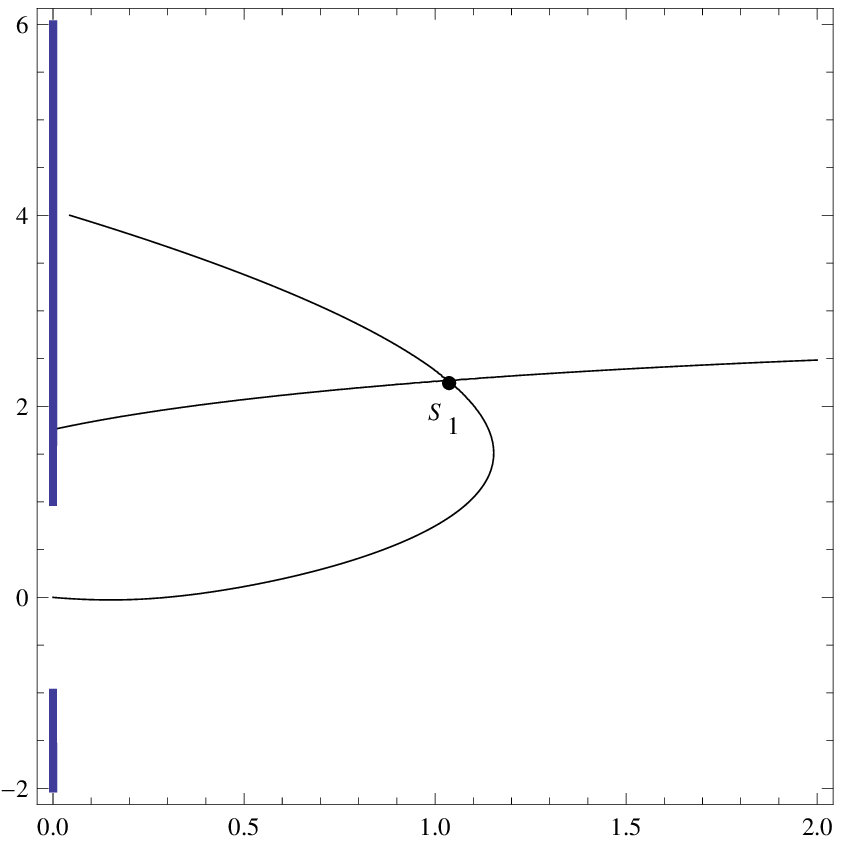}\qquad {\tiny($b$)} \includegraphics[width=0.28\textwidth]{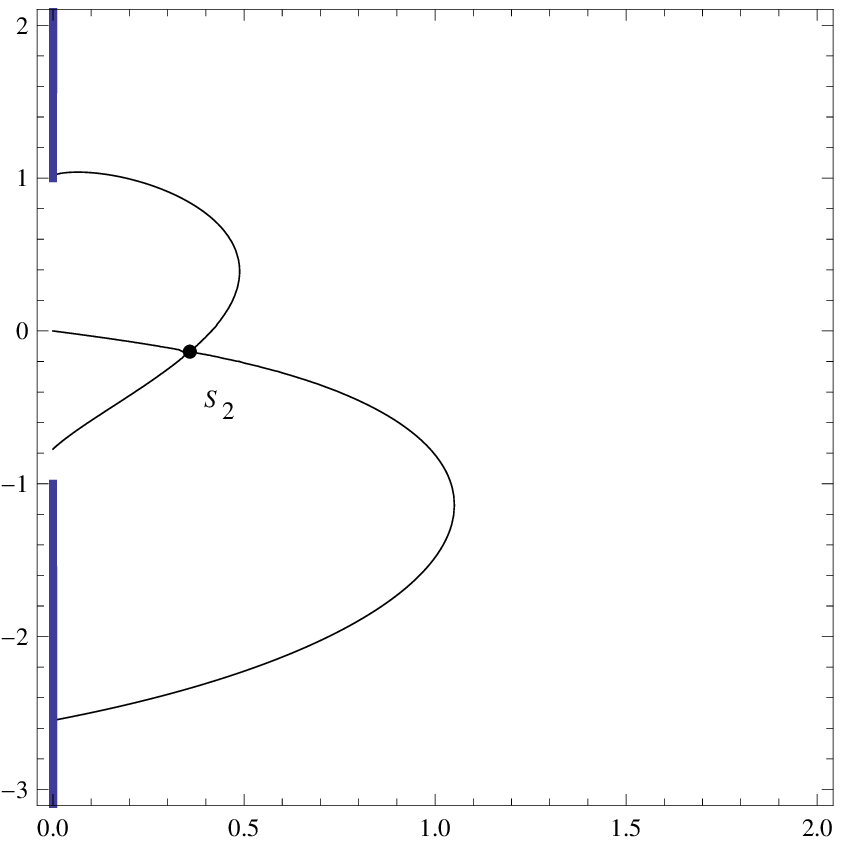} \vspace{.4cm}

{\tiny($c$)}\includegraphics[width=0.28\textwidth]{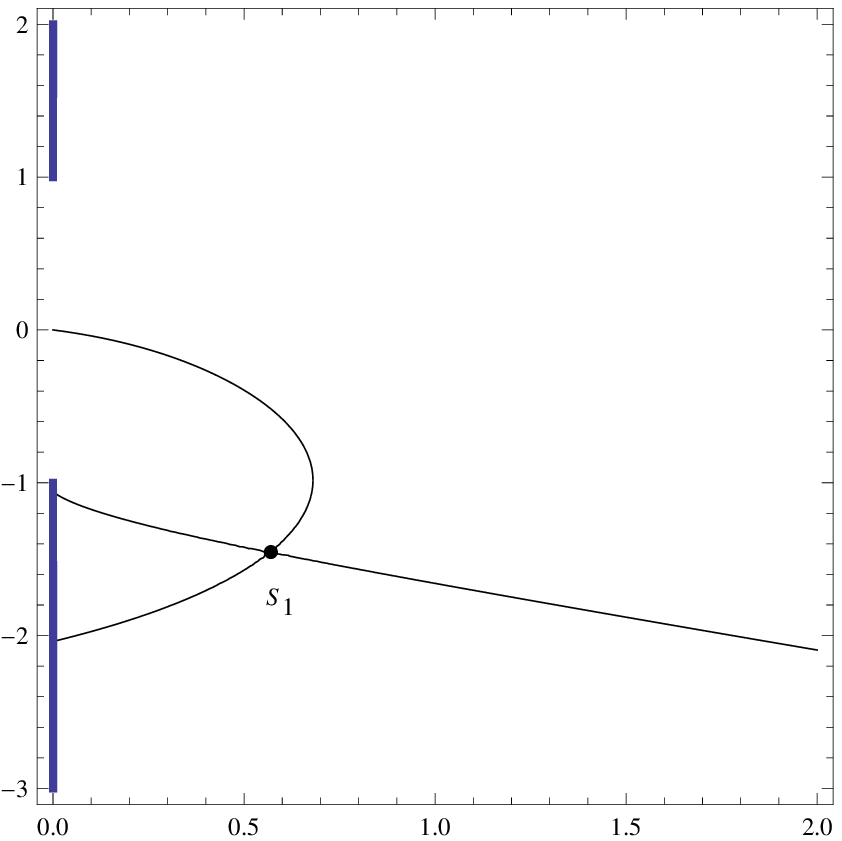}\qquad {\tiny($d$)}\includegraphics[width=0.28\textwidth]{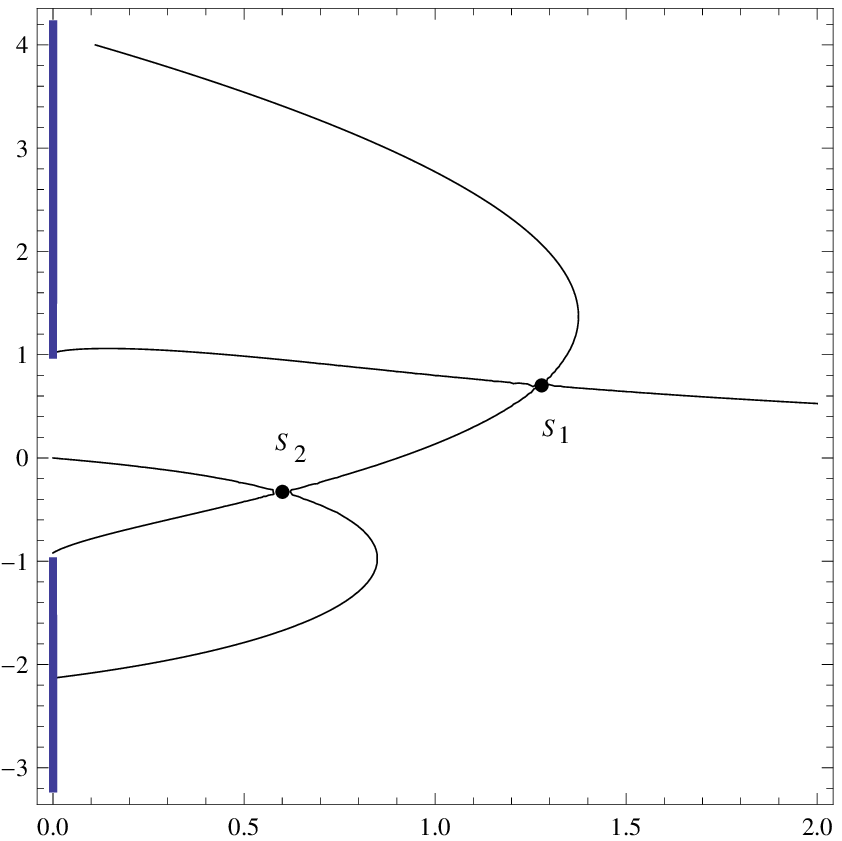}\\
\caption{\small{The steepest paths the through the saddles when $\theta=0.10\pi$: (a) on $PA$ with $q=0.60+0.95307i$, (b) on $PB$ with $q=1.40+0.39447i$, (c) on $PC$ with $q=0.40-0.42914i$ and (d) at $P$ with $q=0.93778+0.18745i$. The heavy lines denote the branch cuts.}}
	\end{center}
\end{figure}
\vspace{0.6cm}

\begin{center}
{\bf 4. \ Numerical results}
\end{center}
\setcounter{section}{4}
\setcounter{equation}{0}
\renewcommand{\theequation}{\arabic{section}.\arabic{equation}}
To verify these assertions, we carry out calculations using (\ref{e22}) and (\ref{e23}) for a series of values of $q\equiv\nu/z$ situated in different domains in Fig.~2. The results are presented in Table 3 which shows the absolute relative error in the computation of ${\bf H}_\nu(z)$. The values of the Bessel functions $J_\nu(z)$ and $Y_\nu(z)$ were evaluated with the in-built codes in {\it Mathematica}. In each case, the asymptotic series on the right-hand sides of (\ref{e22}) and (\ref{e23}) is optimally truncated; that is, at or just before the least term.  
\begin{table}[th]
\caption{\footnotesize{The absolute relative error in the computation of ${\bf H}_\nu(z)$ from (\ref{e22}) and (\ref{e23}) when $z=40e^{i\theta}$. }}
\begin{center}
\begin{tabular}{|l|lc|lc|}
\hline
\mcol{1}{|c|}{} & \mcol{2}{c|}{$\theta=0$} & \mcol{2}{c|}{$\theta=0.10\pi$}\\
\mcol{1}{|c|}{$q=\nu/z$} & \mcol{1}{c}{Error} & \mcol{1}{c|}{Endpoint} & \mcol{1}{c}{Error} & \mcol{1}{c|}{Endpoint}\\
[.1cm]\hline
&&&&\\[-0.3cm]
$0.60$ & $7.764\times 10^{-9}$ & $\infty$ & $9.556\times 10^{-9}$ & $\infty$\\
$1.00$ & $1.041\times 10^{-4}$ & $\pm i$ & $2.751\times 10^{-5}$ & $-i$\\
$1.25$ & $8.835\times 10^{-4}$ & $\pm i$ & $4.830\times 10^{-4}$ & $-i$\\ 
$0.60+0.40i$ & $2.355\times 10^{-6}$ & $\infty$ & $3.280\times 10^{-6}$ & $\infty$\\ 
$1.00+0.60i$ & $2.783\times 10^{-4}$ & $+i$ & $4.136\times 10^{-3}$ & $+i$\\
$1.00-0.30i$ & $7.342\times 10^{-5}$ & $-i$ & $5.000\times 10^{-5}$ & $-i$\\
[.2cm]\hline
\end{tabular}
\end{center}
\end{table}

In \cite[\S 10.43]{W}, Watson claims that (\ref{e23}) and (\ref{e22}) hold for $q\in (0,1)$ and $q\in [1,\infty)$, respectively, when $|\arg\,z|<\fs\pi$. Our calculations have shown that when $\arg\,z\neq 0$ with $q=\nu/z>0$ (that is, when $\nu$ and $z$ have the same phase), the expansion (\ref{e23}) holds for $q\in (0,Q)$ and the expansion (\ref{e22}) holds for $q\in [Q,\infty)$, where $Q\equiv Q(\theta)$ is tabulated in Table 2.

\end{document}